\documentclass[leqno]{article}
\usepackage{amsmath,amsfonts,amsthm}

\newtheorem{prop}{Proposition}[section]
\newtheorem{defn}[prop]{Definition}

\newtheorem{conj}[prop]{Conjecture}

\newtheorem{rem}[prop]{Remark}

\newcommand\Z{\mathbb{Z}}
\newcommand\C{\mathbb{C}}
\newcommand\Q{\mathbb{Q}}
\newcommand\p{\mathbb{P}}

\newcommand\cL{{\mathcal{L}}}
\newcommand\cO{{\mathcal{O}}}
\newcommand\cQ{{\mathcal{Q}}}

\newcommand\cV{{\mathcal{V}}}

\newcommand\M[2]{\overline{M}_{#1,#2}}

\renewcommand\b{\beta}
\newcommand\g{\gamma}
\renewcommand\t{\tau}

\newcommand\hidot{{\raise1pt\hbox{$\scriptscriptstyle\bullet$}}}
\newcommand\lodot{{\raise.3pt\hbox{$\scriptscriptstyle\bullet$}}}

\begin{document}

\title{Virtual Fundamental Classes of Zero Loci}

\author{David A.\ Cox\\
Department of Mathematics and Computer Science\\
Amherst College\\
Amherst, MA 01002-5000\\
{\tt dac@cs.amherst.edu}\\
\and
Sheldon Katz\\
Department of Mathematics\\
Oklahoma State University\\
Stillwater, OK 74078-0613\\
{\tt katz@math.okstate.edu}\\
\and
Yuan-Pin Lee\\
Department of Mathematics\\
University of California at Los Angeles\\
Los Angeles, CA 90095-1555\\
{\tt yplee@math.ucla.edu}}

\date{}

\maketitle

\section{Introduction and a conjecture}
\label{statement}

Let $X$ be a smooth projective variety over $\C$.  For each
nonnegative integer $n$ and homology class $\b\in H_2(X,\Z)$, we are
interested in $n$-pointed genus $0$ stable maps to $X$ of class $\b$.
The concept of a stable map was introduced in \cite{kontsevich}.  We
will for the most part follow the notation of \cite{ck}.  In
particular, such a stable map is
\[
f:(C,p_1,\ldots,p_n)\to X,
\]
where $p_a(C)=0$ and $f_*[C]=\b$.  Families of stable maps can be
described either by a coarse moduli space or by a stack.  Here, we
will use the language of stacks, and we will denote this stack by
$\M0{n}(X,\b)$.  For more about stacks, see \cite{delmum,stacks}.
Roughly speaking, $\M0{n}(X,\b)$ is a functor
\[
\M0{n}(X,\beta) : (\C\hbox{-{\rm Schemes}}) \longrightarrow (\hbox{\rm
Sets}) 
\]
which associates to any scheme $S$ over $\C$ the set of all
isomorphism classes of families of $n$-pointed genus~0 stable maps to
$X$ over $S$ in the class $\b$.  (As with many moduli problems,
$\M0{n}(X,\b)$ is actually a groupoid rather than a functor.  This is
explained in the appendix to \cite{vistoli}.)  One can show that
$\M0{n}(X,\b)$ is a Deligne-Mumford stack \cite{behrendmanin}.  The
associated coarse moduli space was shown to exist in \cite{alexeev}
and was explicitly constructed in \cite{fp}.

It is well known that the stack $\M0{n}(X,\b)$ can have pathological
behavior.  For instance, even its dimension does not remain constant
under deformation of complex structure of $X$.  The problem is due to
obstructions, which arise in a very general setting and can behave 
erratically.

Specifically, the tangent space to $\M0{n}(X,\b)$ at
$f:(C,p_1,\ldots,p_n)\to X$ is
\[
\mathrm{Ext}^1_C\big(f^*\Omega^1_X \to
\Omega^1_C({\textstyle{\sum_{i=1}^n}} p_i),\mathcal{O}_C\big),
\]
while the obstructions lie in
\[
\mathrm{Ext}^2_C\big(f^*\Omega^1_X \to
\Omega^1_C({\textstyle{\sum_{i=1}^n}} p_i),\mathcal{O}_C\big).
\]
The dimensions of these $\mathrm{Ext}^1$s and $\mathrm{Ext}^2$s can vary
as the complex structure of $X$ varies, but the difference 
$\dim(\mathrm{Ext}^1)-\dim(\mathrm{Ext}^2)$
remains constant, equal to
\[
d=\int_\b c_1(X)+n+\dim X-3.
\]
This quantity is called the {\em virtual dimension\/} of $\M0{n}(X,\b)$.

These considerations lead to natural notions of an obstruction theory,
and an associated virtual fundamental
class  \cite{litian,behrfant}
\[
[\M0{n}(X,\b)]^{\rm virt}\in A_d(\M0{n}(X,\b)).
\]
The virtual fundamental class is invariant under deformations of the
complex structure of $X$ in a precise sense.  For example, if $X$ is a
Calabi-Yau threefold and $n=0$, then the virtual dimension is 0, and
we have the Gromov-Witten invariant
\[
N_\b=\mathrm{deg}\,[\M00(X,\b)]^{\mathrm{virt}}\in\Q.
\]
This number is independent of the complex structure of $X$, and is
closely related to the ``number of rational curves'' on $X$ in the
homology class $\b$.  One approach to making sense of the ``number of
curves'' is to introduce the notion of an {\em instanton number\/}
$n_\b$ defined recursively by
\begin{equation}
\label{instanton}
N_\b=\sum_{k\mid\b}\frac{n_{\b/k}}{k^3}.
\end{equation}
The individual terms in~(\ref{instanton}) morally account for the
contribution to $N_\b$ of degree $k$ covers of $n_{\b/k}$ distinct
embedded curves with homology class $\b/k$.  The $n_\b$ are
conjectured to be integers.  However, even if this integrality were
proven, these numbers are not quite the same as the elusive ``number
of rational curves''.\footnote{The $n_\b$ and their generalizations to
higher genus have a direct meaning in physics \cite{gv}.  If these
ideas could be transformed into rigorous mathematics, integrality
would follow immediately.  Some evidence is presented in
\cite{bkl,bp,kkv}.} 
See \cite[Section 7.4.4]{ck} for further discussion.

We are now ready to formulate our conjecture.  For each $k=1,\ldots,
n$ we have the usual evaluation maps $e_k:\M0{n}(X,\b)\to X$, which on
geometric points take $f:(C,p_1,\ldots,p_n)\to X$ to $f(p_k)$.
Furthermore, the universal stable curve over $\M0{n}(X,\b)$ is the map
$\pi_{n+1}:\M0{n+1}(X,\b) \to \M0{n}(X,\b)$ which ignores the last
marked point and contracts any components which have become unstable.
In Section~\ref{relations}, we will also need to use the $n$
tautological sections $s_j$ for $j=1,\ldots n$, where $s_j$ takes a
pointed stable map to the point $p_j$, identifying the fiber of the
universal curve over $f:(C,p_1,\ldots, p_n)\to X$ with $C$.  Further
details of this situation are given in \cite[Section 10.1.1]{ck}.

Now let $V$ be a convex vector bundle on $X$, so that $H^1(f^*V)=0$
for all genus~0 stable maps $f$ to $X$.  We then put
\[
\cV_{\b,n}=(\pi_{n+1})_*e_{n+1}^*V, 
\]
which is a vector bundle on $\M0n(X,\b)$ by the convexity of $V$.  Let
$i:Y\hookrightarrow X$ be the inclusion of the zero locus of a regular
section of $V$.  For any homology class $\g\in H_2(Y,\Z)$, this
induces a natural inclusion $j_\g:\M0{n}(Y,\g)\hookrightarrow
\M0{n}(X,i_*\g)$.

\begin{conj}  
\label{virtconj}
For any $\b\in H_2(X,\Z)$, we have
\begin{equation}
\label{e:virtconj}
\sum_{i_*{\g}=\b}(j_\g)_*[\M0n(Y,\g)]^{\rm virt}=
c_\mathrm{top}(\cV_{\beta,n})\cap [\M0n(X,\b)]^{\rm virt}.
\end{equation}
\end{conj}

Here and throughout the rest of the paper, we use the notation
$c_\mathrm{top}(E)$ to denote the top chern operator of the vector
bundle $E$.

As a special case, this conjecture explains Kontsevich's original
formulation of the Gromov-Witten invariants of the quintic threefold.
Here, $X=\p^4,\ V=\cO(5)$, and $Y$ is a quintic threefold.  Homology
classes on $Y$ and on $\p^4$ can be identified with the degree $d$.
Putting $n=0$, Conjecture~\ref{virtconj} reads
\[
(j_d)_*[\M00(Y,d)]^{\rm virt}=c_{\mathrm{top}}(\cV_{d,0})\cap
[\M00(\p^4,d)]={\rm Euler}(\cV_{d,0}),
\]
since $[\M00(\p^4,d)]^{\rm virt}$ is the usual fundamental class by
the convexity of $\p^4$.  We then have
\begin{equation}
\label{nd}
N_d=\mathrm{deg}\,[\M00(Y,d)]^{\mathrm{virt}}=
\int_{\M00(\p^4,d)}{\rm Euler}(\cV_{d,0}),
\end{equation}
agreeing with Kontsevich's definition~\cite{kontsevich}.  A proof of
the conjecture in this case can be found in \cite[Example
7.1.5.1]{ck}. 

The intuition behind Conjecture \ref{virtconj} is easy to explain.  If
$Y$ is the zero locus of a regular section $s$ of $V$, then (roughly
speaking), we should have the following:
\begin{itemize}
\item The fiber of $\cV_{d,0}$ over $f : C \to
X$ is $H^0(f^*(V))$.  Thus $f \mapsto f^*(s)$ is a section of
$\cV_{d,0}$ over $\M0n(X,\b)$.
\item This section vanishes when $f$ factors through $i : Y
\hookrightarrow X$.  Thus the
zero locus of this section is $\coprod_{i_*\gamma=\b} \M0n(Y,\gamma)$.
\item At the level of fundamental classes, the zero locus of a regular
section of a vector bundle is given by evaluation of its top chern
operator on the fundamental class.
\end{itemize}
In a few cases (including the quintic three discussed above), this
naive reasoning can be made rigorous.  In the general case, the
presence of the virtual fundamental class makes the situation more
complicated, though one expects the conjecture to follow easily once
the general properties of the virtual fundamental class are fully
understood.  We will give a special case of how this works in the next
section.

We should also mention that Conjecture~\ref{virtconj} generalizes the
conjecture stated in \cite[(11.81)]{ck}.  As will we see in
Section~\ref{relations}, this has implications for the Gromov-Witten
invariants of certain nef complete intersections in smooth toric
varieties.

\section{Proof in a special case}
\label{special}

The basic idea of this section is that Conjecture~\ref{virtconj} is
true when $V$ is an \emph{embedding vector bundle} on $X$.  This is
defined as follows.  Let $G(r,k)$ be the Grassmannian of $r$
dimensional quotients of $\C^k$.  

\begin{defn}
\label{positive}
A vector bundle $V$ of rank $r$ on $X$ is an {\bf embedding bundle} if $V$
is generated by global sections and the map $X \to G(r,k)$, $k =
h^0(X,V)$, induced by the exact sequence
\[ 
H^0(X,V)\otimes \cO_X \to V \to 0
\]
is an embedding.
\end{defn}

In particular, when $V$ is an embedding bundle, it is the restriction of
the universal quotient bundle $\cQ$ on $G(r,k)$.  Of course, when $r =
1$, ``embedding'' is equivalent to ``very ample''.  However, when
$r > 1$, ``embedding'' differs from the notion of ``positive'' or
``ample'' vector bundle defined in \cite{griffiths} or
\cite{hartshorne}.

We can now prove our conjecture in the case of an embedding bundle.

\begin{prop}
\label{specialcase.2}
Conjecture~\ref{virtconj} is true when the vector bundle $V$ on $X$ is
embedding and $Y \subset X$ is the zero locus of a generic section
of $V$.
\end{prop}

\begin{proof}
We will use an argument in \cite{gathmann} together with standard
facts about refined Gysin maps from \cite{fulton}. The rough idea of
the proof is that the conjecture is obviously true for a Grassmannian
$G(r,k)$, and then the results of \cite{gathmann} and \cite{fulton}
will show that it remains true when we pull back from $G(r,k)$ to $X$.

To see how this works in detail, we first recall some notation.
Suppose we have a cartesian diagram of Deligne-Mumford stacks:
\begin{equation}
\label{diagram1.2}
\begin{matrix}
M & \xrightarrow{\,j\,} & M_2 \\[2pt]
\downarrow && \phantom{\scriptstyle f} \downarrow {\scriptstyle f} \\
M_1 & \xrightarrow{\,i\,} & S
\end{matrix}
\end{equation}
where $M_1 \to S$ is a regular embedding.  In this situation, we have
the refined Gysin map $i^! : A_*(M_2) \to A_*(M)$ defined in
\cite[Chapter 6]{fulton}, where we use \cite{vistoli} to extend from
schemes to stacks.  

Now assume that there is also a vector bundle $\cV$ on $S$ such that
$i : M_1 \subset S$ is the inclusion of the zero locus of a regular
section of $\cV$.  By \cite[Ex.\ 6.3.4]{fulton}, it follows that
\begin{equation}
\label{fultonprop}
j_*i^!(\g) = c_\mathrm{top}(f^*\cV)\cap\g\
\end{equation}
for all $\g \in A_*(M_2)$, where $j$ and $f$ are from
\eqref{diagram1.2}.

To apply this to our situation, suppose that $V$ is an embedding
bundle on $X$ and $Y \subset X$ is the zero locus of a section $s$ of
$V$.  By Definition~\ref{positive}, we can assume that $X$ is embedded
in the Grassmannian $G(r,k)$.  Furthermore, $s$ induces a section
$s_\cQ$ of the universal quotient bundle $\cQ$ on $G(r,k)$ via the
tautological quotient mapping
\[
H^0(X,V)\otimes\cO_{G(r,k)}\to \cQ.
\]
Let $G \subset G(r,k)$ be the zero locus of $s_\cQ$.  It follows from
this description that $s_\cQ$ is a regular section of $\cQ$.  This
description also implies that $G \simeq G(r,H^0(X,V)/\C\cdot s) \simeq
G(r,k-1)$ and that $Y = X\cap G$.

Now fix $\beta \in H_2(X)$ and suppose that $\beta$ maps to $d \in
H_2(G(r,k)) \simeq H_2(G) \simeq \Z$.  Then we have a cartesian diagram:
\begin{equation}
\label{cartesianYX}
\begin{matrix}
\coprod_{i_*\g = \beta} \M0n(Y,\g)\! & \xrightarrow{\,j\,} &
\M0n(X,\beta) \\[2pt] 
\downarrow\! && \phantom{\scriptstyle f}\downarrow {\scriptstyle f}
\\[2pt] 
\M0n(G,d)\! & \xrightarrow{\,i\,} & \M0n(G(r,k),d),
\end{matrix}
\end{equation}
where $i$ and $f$ are the natural inclusions.

Since $G(r,k)$ and $G \simeq G(r,k-1)$ are homogeneous spaces, it
follows that $\M0n(G,d) \to \M0n(G(r,k),d)$ is a regular embedding of
smooth stacks.  Applying the above construction, we get the class
\[
i^!([\M0n(X,\beta)]^{\rm virt}) \in \sum_{i_*\gamma =\beta}
A_*(\M0n(Y,\g)). 
\]
Using the argument of \cite[Lemma 4.2]{gathmann}, we see that
\begin{equation}
\label{gathmann}
i^!([\M0n(X,\beta)]^{\rm virt}) = \sum_{i_*\gamma =\beta}
[\M0n(Y,\g)]^{\rm virt}. 
\end{equation}
Gathmann's lemma uses $\p^N$ and $H \simeq \p^{N-1}$ rather than
$G(r,k)$ and $G \simeq G(r,k-1)$, but the proof still applies since
$G(r,k)$ and $G$ are convex.  Also, while the statement of Lemma 4.2
in \cite{gathmann} is different from \eqref{gathmann}, the final
sentence of his proof shows that \eqref{gathmann} follows from his
argument.

In Section~\ref{statement}, we discussed how $V$ induces the bundle
$\cV_{\b,n}$ on $\M0n(X,\b)$.  In a similar way, the universal
quotient bundle $\cQ$ induces a bundle $\cV_{d,n}=(\pi_{n+1})_*
e_{n+1}^*\cQ$ on $\M0n(G(r,k),d)$. The map $i : \M0n(G,d)
\to \M0n(G(r,k),d)$ is the zero locus of a section of $\cV_{d,n}$ by
\cite[Sect.~2.1]{pandbour}.  Using \eqref{fultonprop}, it follows that
\[
j_*i^!(\g) = c_\mathrm{top}(f^*\cV_{d,n})\cap\g
\]
for all $\g \in A_*(\M0n(X,\beta))$, where $j$ and $f$ are from
\eqref{cartesianYX}.  Since $f^*\cV_{d,n} = \cV_{\b,n}$, we obtain
\[
j_*i^!\big([\M0n(X,\beta)]^{\rm virt}\big) =
c_\mathrm{top}(\cV_{\b,n})\cap[\M0n(X,\b)]^\mathrm{virt}.
\]
Combining this with \eqref{gathmann}, the proposition follows.
\end{proof}

\begin{rem} {\rm The obvious embedding
\[
G(r_1,k_1) \times G(r_2,k_2) \hookrightarrow G(r_1+r_2,k_1+k_2)
\]
shows that a direct sum of embedding bundles is embedding.  In
particular, a direct sum of very ample line bundles is an embedding
bundle.  Thus Conjecture~\ref{virtconj} holds when $Y$ is a generic
complete intersection of very ample hypersurfaces in $X$.

In some versions of the mirror theorem, the ambient space is a
smooth toric variety.  Here, a line bundle is very ample if and only
if it is ample, so that Conjecture~\ref{virtconj} holds when $Y$ is a
generic complete intersection of ample hypersurfaces in a smooth toric
variety $X$.

In Remark~4.5 of \cite{gathmann}, Gathmann says that he expects
that his results should hold under the weaker hypothesis that $V$ is
generated by global sections.  (Gathmann only considers
line bundles, but if his argument extends to line bundles generated by
global sections, then it should also work for vector bundles generated
by global sections.)

For a smooth (or simplicial) toric variety $X$, one can easily
show that a line bundle is convex if and only if it is generated by
global sections (this follows from Reid's description
\cite[Prop.~1.6]{reid} of the Mori cone of $X$).  Thus, if Gathmann's
result can be extended to vector bundles generated by global sections,
then it would follow that Conjecture~\ref{virtconj} would hold when
the vector bundle $V$ is a direct sum of convex line bundles on a
smooth toric variety $X$.}
\end{rem}

\section{Relation to mirror theorems}
\label{relations}

This section will discuss the relationship between
Conjecture~\ref{virtconj} and various approaches to the mirror
theorem.

This basic idea is that computing Gromov--Witten invariants of
hypersurfaces $Y\subset X$ (or more generally, the zero locus of a
section of a vector bundle $V$ on $X$) is a two-stage process:
\begin{itemize}
\item First, one must relate the Gromov--Witten invariants of $Y$ to
invariants defined using the ambient space $X$ and the vector bundle
$V$.  Conjecture~\ref{virtconj} shows how to do this.
\item Second, one must compute the invariants defined using $X$ and
$V$.  This is what the various \emph{mirror theorems} in the
literature do.
\end{itemize}
We will briefly discuss three approaches to the mirror theorem
and describe the role of equation \eqref{e:virtconj} from
Conejecture~\ref{virtconj} in each of these (related) approaches.

Givental's approach~\cite{givental1} was largely based on the fact
that the quantum cohomology of $X$ can be described by a quantum
$\mathcal{D}$-module generated by a single formal function, which we
denote by $J_X$.  We recall the definition, following
\cite[Chapters~10 and 11]{ck}.

Let $\{T_0,\ldots T_m\}$ be a basis for $H^*(X)$ 
with $T_0=1$, and $\{T_1,\ldots,
T_r\}$ a basis for $H^2(X)$.  Let
$\{T^a\}$ be the dual basis under the intersection pairing.  Introduce
variables $t_i$ and put $\delta=\sum_{i=1}^rt_iT_i$.  

We recall the definition of the gravitational correlators.  The
universal curve $\M02(X,\b) \to \M01(X,\b)$ has the
section $s_1$ as described in Section~\ref{statement}.
We put $\cL=s_1^*(\omega)$, where $\omega$ is the relative dualizing
sheaf of $\M02(X,\b) \to \M01(X,\b)$.  Then the 1-point genus 0
gravitational correlators are defined by the equation
\begin{equation}
\label{grcor}
\langle \tau_n T_a\rangle_{0,\b}=\int_{\M01(X,\b)}c_1(\cL)^n\cup
e_1^*(T_a) \cap [\M01(X,\b)]^{\rm virt}.
\end{equation}

Let $q^\b$ be formal symbols satisfying $q^\b\cdot q^{\b'}=
q^{\b+\b'}$ and let $\hbar$ be a formal parameter.  Then we define the
formal function $J_X$ by the equation
\begin{equation}
\label{jdef}
J_X= e^{(t_0+\delta)/\hbar}\left(
1+\sum_{\b\ne0}\sum_{a=0}^m\sum_{n=0}^\infty
\hbar^{-(n+2)}\langle\t_nT_a\rangle_{0,\b}\,T^aq^\b\right).
\end{equation}

The goal is to relate $J_Y$ to a variant of $J_X$ which takes the
bundle $V$ into account.  We denote the kernel of the natural
``evaluation map'' $\cV_{\b,1}\to e_1^*V$ by $\cV'_{\b,1}$.  We then
put
\begin{equation}
\label{jv}
J_{\cV}=e^{(t_0+\delta)/\hbar}\mathrm{Euler}(\mathcal{V})\bigg(
1+\sum_{\b\ne0}q^\b\, e_{1!}\Big(
\frac{\mathrm{Euler}(\mathcal{V}'_{\b,1})}{\hbar(\hbar-c)}
\Big)\bigg),
\end{equation}
where the expression $1/(\hbar-c)$ is understood to be expanded as
\[
\frac1{\hbar-c}=\sum_{n=0}^\infty c^n\hbar^{-(n+1)}
=\sum_{n=0}^\infty c_1(\cL)^n\hbar^{-(n+1)}.
\]
It was shown in \cite{ck} Section~11.2 that when $X = \p^n$,
Conjecture~\ref{virtconj} implies that
\begin{equation}
\label{ijyjv}
i_*J_Y=J_{\cV},
\end{equation}
relating the Gromov-Witten invariants on $Y$ to invariants defined
using $\p^n$ and $V$, as desired.

Givental's approach to the mirror theorem computes $J_\cV$ explicitly
in terms of hypergeometric functions.  This is described in
\cite{givental2} (see the references therein) and Theorem~11.2.16 of
\cite{ck}, and is extended by the Quantum Hyperplane Section Principle
of \cite{kim}.  Once one has $J_\cV$, one needs a formula such as
\eqref{ijyjv} in order to compute the Gromov-Witten invariants of $Y$.

For convex bundles on projective space, a result analogous to
\eqref{ijyjv} is proved in \cite{ckyz}, with an extension to direct
sums of convex/concave line bundles on projective space in
\cite{elezi}.  Hence \eqref{ijyjv} and its analogues provide a bridge
between Gromov-Witten invariants on $Y$ and hypergeometric functions.
Once Conjecture~\ref{virtconj} is proved, we expect that similar
formulas should hold whenever $Y$ is the zero locus of a regular
section of a convex vector bundle $V$ on $X$.  This will be useful
when generalizing the results of \cite{ckyz,elezi} to other ambient
spaces.

Another approach to the mirror conjecture is the ``Mirror Principle''
proposed by Lian, Liu and Yau \cite{lly1,lly2,lly3,lly4}.  In the
third of these papers, the basic object of interest is the integral
\begin{equation}
\label{lly}
K_\beta = \int_{[\overline{M}_{0,0}(X,\beta)]^{\rm virt}} b(\cV_{\beta,0}),
\end{equation}
where:
\begin{itemize}
\item $X$ is a projective manifold and $\b \in H_2(X)$.
\item $V$ is the concavex vector bundle on $X$.  (A \emph{concavex
bundle} is a direct sum of a convex bundle and a concave bundle.)
\item When $V$ is convex, $\cV_{\beta,0}$ is the vector bundle
$\pi_*e_{1}^*V$ on $\M00(X,\beta)$ as defined in
Section~\ref{statement}.  When $V$ is concave, $\cV_{\beta,0}=
R^1\pi_*e_{1}^*V$.  
\item $b$ is a multiplicative characteristic class, such as
the Euler class.
\end{itemize}

The paper \cite{lly3} discusses the properties of the generating
function which has the $K_\b$ as coefficients.  Explicit formulas for
this generating function are given when $X$ is a balloon manifold and
$b$ is either the Euler class or the Chern polynomial \cite[section
4]{lly3}.

It follows that there are many situations where the $K_\b$ can be
computed.  This raises the question of interpreting these numbers in
terms of Gromov-Witten invariants.  Let us explain how this works when
$V$ is a convex bundle on $X$ and $Y \subset X$ is the zero locus of a
regular section of $V$.

In this situation, suppose that the multiplicative characteristic
class $b$ is the Euler class and assume also that
Conjecture~\ref{virtconj} holds for $V$ and $Y \subset X$.  Then the
integral $K_\b$ defined in \eqref{lly} is a sum of Gromov-Witten
invariants of $Y$ as follows:
\[
K_\b = \int_{[\M00(X,\beta)]^{\mathrm{virt}}}\!\!
\mathrm{Euler}(\cV_{\b,0}) = \sum_{i_*\gamma = \beta}
\int_{[\M00(Y,\gamma)]^{\mathrm{virt}}}\!\!1 = \sum_{i_*\gamma =
\beta} \langle I_{0,0,\gamma}\rangle, 
\]
where we are using the notation of \cite[Chapter 7]{ck} for the
Gromov-Witten invariants $\langle I_{0,0,\gamma}\rangle$ of $Y$.
Notice how this generalizes \eqref{nd}.
\smallskip

\begin{rem}
{\rm In the literature, one finds two
algebro-geometric definitions of the virtual fundamental class, one
due to Behrend and Fantechi \cite{behrfant} and the other due to Li
and Tian \cite{litian}.  In Sections~\ref{statement} and
\ref{special}, we used the Behrend-Fantechi definition of virtual
fundamental class.  In particular, the argument of Gathmann used in
the proof of Theorem~\ref{specialcase.2} uses the definition of
$[\M0n(X,\beta)]^{\mathrm{virt}}$ given in \cite{behrfant}.
Lian-Liu-Yau, on the other hand, use the Li-Tian virtual fundamental
class.  However, they never use the explicit construction of Li and
Tian.  What is needed in their papers are the functorial properties of
virtual fundamental classes, which work for Behrend-Fantechi classes
as well.

It is expected that the virtual fundamental classes defined by
Behrend-Fantechi and Li-Tian are equal, but to the best of our
knowledge, the details of this argument have not been written down.}
\end{rem}

There is another proof of the mirror conjecture in \cite{AB} for the
case of complete intersections in $\mathbb{P}^r$.  His approach also
uses the $J$-function of the zero locus $Y$ of a section of $V$ in $X$
via Conjecture~\ref{virtconj}, and relates the $J$-function to
hypergeometric functions.
Let $X=\p^r$ and $V=\mathcal{O}_{\p^r}(l)$ for simplicity. 
Recall that there is a \emph{birational} morphism 
\[
 \varphi: \overline{M}_{0,0}(\p^r \times \p^1, (d,1)) \to N_d
\]
where $N_d =\p^{(r+1)d+r}$ (\cite{ck} Section~11.1.2). In this
approach one first observes that there is a $\mathbb{C}^*$-action on
$\overline{M}_{0,0}(\p^r \times \p^1, (d,1))$ induced from the
$\mathbb{C}^*$-action on the second factor $\p^1$ and that there are
fixed point components which can be identified with
$\overline{M}_{0,1}(\p^r,d)$.  Let
\[ 
 i:\overline{M}_{0,1}(\p^r,d) \hookrightarrow 
  \overline{M}_{0,0}(\p^r \times \p^1, (d,1))
\]
be the inclusion of such a component. Consider the vector bundle
$\mathcal{V}_d=\pi_{1*}e_1^*(V)$ on 
$\overline{M}_{0,0}(\p^r \times \p^1, (d,1))$, where
\[
 \pi_1: \overline{M}_{0,1}(\p^r \times \p^1, (d,1)) \to 
  \overline{M}_{0,0}(\p^r \times \p^1, (d,1))
\]
is the universal curve and
\[
 e_1: \overline{M}_{0,1}(\p^r \times \p^1, (d,1)) \to \p^r
\]
is the evaluation morphism. It is easy to see that 
\[
 c_{\mathrm{top}}(\mathcal{V}_{d,1}) =i^*
 (c_{\mathrm{top}}(\mathcal{V}_{d})). 
\]
The key point is that on the complement of a boundary divisor, one has
an equivariant isomorphism
\[
 \mathcal{V}_{d} \simeq H^0(\p^1,\mathcal{O}(dl)) \otimes
 \varphi^*\bigl(\mathcal{O}_{N_d}(l) \bigr),
\]
which then gives
\[
 c_{\mathrm{top}}(\mathcal{V}_{d})= \varphi^* H_d + \text{boundary
 terms}, 
\]
where $H_d = c_{\mathrm{top}}(H^0(\p^1,\mathcal{O}_{\p^1}(dl)) \otimes
\mathcal{O}_{N_d}(l))$ is an explicit class on the projective space
$N_d$ and the boundary terms are supported on the exceptional divisors
of $\varphi$.  It is then shown that the class $H_d$ gives the
hypergeometric series and that the boundary terms are responsible for
the mirror transformations, which are certain changes of variables
(see \cite[Chapter~11]{ck}).

One novelty of this approach is that, unlike the previous two
approaches, it does not use the torus action on $\p^r$.  It therefore
opens the way to prove mirror theorems for other projective varieties
$X$.  This was recently realized in \cite{YL}, where equation
\eqref{e:virtconj} (in the special case proved in
Proposition~\ref{specialcase.2}) again serves as a starting point of
the proof.

\section{Acknowledgements}

The authors would like to thank Robin Hartshorne, Rahul Pandharipande
and Andrew Sommese for helpful comments regarding terminology.  The
research of the second author was supported in part by NSA grant
MDA904-98-1-0009.  We are also grateful to Emma Previato for
organizing the AMS Special Session in which a preliminary version of
this paper was presented.  Finally, we would like to thank the referee
for useful suggestions.


\begin{thebibliography}{}

\bibitem[Alexeev]{alexeev} V.\ Alexeev, {\em Moduli spaces
$M_{g,n}(W)$ for surfaces\/}, in {\sl Higher-Dim\-ensional Complex
Varieties $($Trento, 1994\/$)$}, de Gruyter, Berlin, 1996, 1--22,
alg-geom/9410003.

\bibitem[BF]{behrfant} K.\ Behrend and B.\ Fantechi, {\em The
intrinsic normal cone\/}, Invent.\ Math.\ {\bf 128} (1997), 45--88,
alg-geom/9601010.

\bibitem[BEFFGK]{stacks} K.\ Behrend, D.\ Edidin, B.\ Fantechi, W.\
Fulton, L.\ G\"ottsche, and A.\ Kresch.  {\sl Introduction to Stacks},
in preparation.

\bibitem[BM]{behrendmanin} K.\ Behrend and Yu.\ Manin, {\em Stacks of
stable maps and Gromov-Witten invariants\/}, Duke J.\ Math. {\bf 85}
(1996), 1--60, alg-geom/9506023.

\bibitem[Bertram]{AB} A.~Bertram, \emph{Another way to enumerate
rational curves with torus actions}, math.AG/ 9905159.

\bibitem[BKL]{bkl} J.~Bryan, S.~Katz, and N.~Conan Leung, {\em
Multiple covers and the integrality conjecture for rational curves in
Calabi-Yau threefolds\/}, math.AG/9911056.

\bibitem[BP]{bp} J.~Bryan and R.~Pandharipande, {\em BPS states of
curves in Calabi-Yau 3-folds\/}, math.AG/ 0009025.

\bibitem[CKYZ]{ckyz} T.-M.\ Chiang, A.\ Klemm, S.-T.\ Yau, and E.\
Zaslow, {\em Local Mirror Symmetry: Calculations and
Interpretations\/}, hep-th/9903053.

\bibitem[CK]{ck} D.\ A.\ Cox and S.\ Katz, {\sl Mirror Symmetry and
Algebraic Geometry\/}, Mathematical Surveys and Monographs {\bf 68},
AMS, Providence, RI, 1999.

\bibitem[DM]{delmum} P.\ Deligne and D.\ Mumford, {\em The
irreducibility of the space of curves of a given genus\/}, Publ.\
Math.\ IHES {\bf 36} (1969), 75--110.

\bibitem[Elezi]{elezi} A.\ Elezi, {\em Mirror symmetry for concavex
vector bundles on projective spaces\/}, math.AG/ 0004157.

\bibitem[Fulton]{fulton} W.\ Fulton, {\sl Intersection Theory\/},
Springer-Verlag, New York-Berlin-Heidelberg, 1984.

\bibitem[FP]{fp} W.\ Fulton and R.\ Pandharipande, {\em Notes on
stable maps and quantum cohomology\/}, in {\sl Algebraic
Geometry---Santa Cruz 1995\/}, Proceedings of Symposia in Pure
Mathematics {\bf 62} (Part 2), AMS, Providence, RI, 1997, 45--96,
alg-geom/9608011.

\bibitem[Gathmann]{gathmann} A.\ Gathmann, {\em Absolute and relative
Gromov-Witten invariants of very ample hypersurfaces\/}, math.AG/9908054.

\bibitem[Givental1]{givental1} A.~Givental, \emph{Equivariant
Gromov-Witten invariants}, IMRN \textbf{13} (1996) 613-663,
alg-geom/9603021.

\bibitem[Givental2]{givental2} A.~Givental, \emph{Elliptic
Gromov-Witten invariants and the generalized mirror conjecture}, in
{\sl Integrable systems and algebraic geometry $($Kobe/Kyoto,
1997\/$)$}, World Sci. Publishing, River Edge, NJ, 1998, 107--155,
math.AG/9803053.

\bibitem[GV]{gv} R.\ Gopakumar and C.\ Vafa, {\em M-theory and
Topological Strings-II\/}, hep-th/9812127.

\bibitem[Griffiths]{griffiths} P.\ Griffiths, {\em Hermitian
differential geometry, Chern classes, and positive vector bundles\/},
in {\sl Global Analysis\/}, Princeton University Press, Princeton, NJ,
1969, 185--251.

\bibitem[Hartshorne]{hartshorne} R.\ Hartshorne, {\em Ample vector
bundles\/}, Publ.\ Math.\ IHES {\bf 29} (1966), 63--94.

\bibitem[KKV]{kkv} S.\ Katz, A.\ Klemm, and C.\ Vafa, {\em M-Theory,
Topological Strings and Spinning Black Holes\/}, 
hep-th/9910181.

\bibitem[Kim]{kim} B.~Kim, \emph{Quantum hyperplane section theorem
for homogeneous spaces}, Acta Math. \textbf{183} (1999), no. 1,
71--99, alg-geom/9712008.

\bibitem[Kontsevich]{kontsevich} M.~Kontsevich, \emph{Enumeration of
rational curves via torus actions}, in {\sl The moduli space of curves
$($Texel Island, 1994\/$)$}, Progress in Math.\ \textbf{129}, Birkh\"auser,
Boston-Basel, 1995, 335--368, hep-th/9405035.

\bibitem[Lee]{YL} Y.-P.~Lee, {\em Quantum Lefschetz hyperplane
theorem}, math.AG/0003128.

\bibitem[LT]{litian} J.\ Li and G.\ Tian, {\em Virtual moduli cycles
and Gromov-Witten invariants of algebraic varieties\/}, J.\ Amer.\
Math.\ Soc.\ {\bf 11} (1998), 119--174, alg-geom/9602007.

\bibitem[LLY1]{lly1} B.\ Lian, K.\ Liu and S.-T.\ Yau, {\em Mirror
Principle I\/}, Asian J.\ Math.\ {\bf 1} (1997), 792--763,
alg-geom/9712011.

\bibitem[LLY2]{lly2} B.\ Lian, K.\ Liu and S.-T.\ Yau, {\em Mirror
Principle II\/}, Asian J.\ Math.\ {\bf 3} (1999), 109--146,
math.AG/9905006. 

\bibitem[LLY3]{lly3} B.\ Lian, K.\ Liu and S.-T.\ Yau, {\em Mirror
Principle III\/}, math.AG/9912038.

\bibitem[LLY4]{lly4} B.\ Lian, K.\ Liu and S.-T.\ Yau, {\em Mirror
Principle IV\/}, math.AG/0007104.

\bibitem[Pandharipande]{pandbour} R.\ Pandharipande, {\em Rational
curves on hypersurfaces $($after A.~Givental\/$)$}, S\'eminaire
Bourbaki, Exp.\ 848, Vol.\ 1997/98, Ast\'erisque {\bf 252} (1998),
307--340, math.AG/9806133.

\bibitem[Reid]{reid} M.\ Reid, {\em Decomposition of toric
morphisms\/}, in {\sl Arithmetic and Geometry\/}, Vol.~II, Progress
in Math.\ {\bf 36}, Birkh\"auser, Boston-Basel, 1983, 395--418.

\bibitem[Vistoli]{vistoli} A.\ Vistoli, {\em Intersection theory on
algebraic stacks and their moduli spaces\/}, Invent.\ Math.\ {\bf 97}
(1989), 613--670.

\end{thebibliography}
\end{document}